\documentclass[12pt]{amsart}
\usepackage[T1]{fontenc}
\usepackage{amssymb}
\usepackage{bm}
\newtheorem{thm}{Theorem}%[section]
\newtheorem{lem}{Lemma}[section]
\newtheorem{prop}{Proposition}[section]

\newtheorem{conj}{Conjecture}
\newtheorem{defi}{Definition}[section]
\newtheorem{ex}{Example}
\newtheorem{rem}{Remark}[section]
\newtheorem{prob}{Problem}%[section]
\usepackage{url}
\usepackage{ulem}
\usepackage[usenames]{color}
\usepackage[utf8]{inputenc}
\numberwithin{equation}{section}

\newcommand{\N}{\mathbb{N}}

\DeclareMathOperator{\lcm}{lcm}

\begin{document}

\title[Number of solutions to $a^x+b^y=c^z$]
{General sharp bounds for the number of solutions to purely exponential equations \\with three terms}

\author{Maohua Le}
\address{Maohua Le
\hfill\break\indent Institute of Mathematics, Lingnan Normal College,
\hfill\break\indent Zhanjiang, Guangdong, 524048.
\hfill\break\indent China
}
\email{lemaohua2008@163.com
}
\author{Takafumi Miyazaki}
\address{Takafumi Miyazaki
\hfill\break\indent Gunma University, Division of Pure and Applied Science,
\hfill\break\indent Graduate School of Science and Technology
\hfill\break\indent Tenjin-cho 1-5-1, Kiryu 376-8515.
\hfill\break\indent Japan
}
\email{tmiyazaki@gunma-u.ac.jp}
\thanks{The second author is supported by JSPS KAKENHI (No. 24K06642).}
%\today
\maketitle

%\markleft{M.-H. Le \& T. Miyazaki}
\markleft{Maohua Le \& Takafumi Miyazaki}
\markright{Number of solutions to $a^x+b^y=c^z$}

%\vspace{-0.3cm}\centerline
%{\sl\footnotesize Dedicated to}

\vspace{-0.0cm}
\begin{abstract}%\bf 
It is conjectured that for any fixed relatively prime positive integers $a,b$ and $c$ all greater than 1 there is at most one solution to the equation $a^x+b^y=c^z$ in positive integers $x,y$ and $z$, except for specific cases.
In this paper, we prove that for any fixed $c$ there is at most one solution to the equation, except for only finitely many pairs of $a$ and $b.$ 
This is regarded as a 3-variable generalization of the result of Miyazaki and Pink [T. Miyazaki and I. Pink, Number of solutions to a special type of unit equations in two unknowns, I\hspace{-1.2pt}I\hspace{-1.2pt}I, arXiv:2403.20037 (accepted for publication in Math. Proc. Cambridge Philos. Soc.)] which asserts that for any fixed positive integer $a$ there are only finitely many pairs of coprime positive integers $b$ and $c$ with $b>1$ such that the Pillai's type equation $a^x-b^y=c$ has more than one solution in positive integers $x$ and $y$.
The proof of our result is based on a certain $p$-adic idea of Miyazaki and Pink and relies on many deep theorems on the theory of Diophantine approximation, and it also includes the complete description of solutions to some interesting system of simultaneous polynomial-exponential equations. 
We also discuss how effectively exceptional pairs of $a$ and $b$ on our result for each $c$ can be determined.
\end{abstract}

%%%%%%%%%%%%%%%%%%%%%%%%%%%
\section{Introduction}\label{sec-intro}
%%%%%%%%%%%%%%%%%%%%%%%%%%%

An important problem in Diophantine number theory is to establish sharp upper bounds for the number of solutions to unit equations (cf.~\cite{EvGy}).
In this paper we deal with a class of $S$-unit equations in two unknowns over the rationals which are closely related to Pillai's conjecture (cf.~\cite{Wa}) and generalized Fermat conjecture (cf.~\cite{Beu,BenMihSik}).

In a recent decade, there have been important progresses on the problem of establishing general sharp upper bounds for the number of solutions to the following exponential Diophantine equation:
\begin{equation} \label{abc}
a^x+b^y=c^z
\end{equation}
in positive integers $x,y$ and $z$, where $a,b$ and $c$ are fixed relatively prime positive integers all greater than 1.
The above equation is a kind of basic and important purely exponential equations, and the study of it has a long and rich history with the pioneer works of Mahler \cite{Ma}, Gel'fond \cite{Ge} and Siepi\'nski \cite{Sie} (see also \cite[Ch.\,1]{ShTi}, \cite[Sec.\,D10]{Gu}, \cite[Ch.\,4--6]{EvGy} or \cite{MiyPi,MiyPi2,MiyPi3}).
We are interested in finding a general and sharp upper bound for % $N=N(a,b,c)$, 
the number of solutions $(x,y,z)$ to \eqref{abc}.

Among the recent progresses referred to before, Miyazaki and Pink \cite{MiyPi} finally proved the following definitive result:

\begin{prop}%[Theorem 1 of \cite{MiyPi}] \label{atmost2}
There are in general at most two solutions to equation $\eqref{abc},$ except for $(a,b,c)=(3,5,2)$ or $(5,3,2)$ each of which gives exactly three solutions.
%$N(a,b,c) \le 2,$ except for $N(3,5,2)=N(5,3,2)=3.$% each of which gives exactly three solutions.
\end{prop}

This result is sharp in the sense that there are infinitely many cases in which equation \eqref{abc} has exactly two solutions %with $N=2$ 
(see \eqref{ex-excep-abc} below).

On the other hand, something rather stronger seems to be true, and a formulation in that direction was put forward by Scott and Styer \cite{ScSt_PMD_2016}, as follows:

\begin{conj}%[\cite{ScSt_PMD_2016}]
\label{atmost1}
Assume that none of $a,b$ and $c$ is a perfect power.
Then there is at most one solution to equation $\eqref{abc},$ except when $(a,b,c)$ or $(b,a,c)$ belongs to the following set\,$:$
\begin{align}\label{excep-set}
\{ \,&(3,5,2),(3,13,2),(2,5,3),(2,7,3),\\
&(2,3,11),(3,10,13),(2,3,35),(2,89,91),\nonumber\\
&(2,5,133),(2,3,259),(3,13,2200),(2,91,8283),\nonumber\\
&(2,2^r-1,2^r+1)\, \},\nonumber
\end{align}
where $r$ is any positive integer with $r=2$ or $r \ge 4.$
\end{conj}

This conjecture is regarded as a 3-variable generalization of the conjecture of Bennett \cite[Conjecture 1.2]{Ben_cjm_01} on Pillai's type equation, and it is one of the utmost importance throughout the study of purely exponential Diophantine equations.
It is known that the solutions to equation \eqref{abc} corresponding to the cases in \eqref{excep-set} are found from the following set of equalities:
\begin{gather}
3^{}+5^{}=2^{3}, \ 3^{3}+5^{}=2^{5}, \ 3^{}+5^{3}=2^{7}; \nonumber\\
3^{}+13^{}=2^{4}, \ 3^{5}+13^{}=2^{8}; \nonumber\\
%1^{}+2^{}=3^{}, \ 1^{}+2^{3}=3^{2}; \nonumber\\
2^{2}+5^{}=3^{2}, \ 2^{}+5^{2}=3^{3}; \nonumber\\
2^{}+7^{}=3^{2}, \ 2^{5}+7^{2}=3^{4}; \nonumber\\
2^{3}+3^{}=11^{}, \ 2^{}+3^{2}=11^{}; \nonumber\\ %<== 1+2=3, 1+2^3=3^2
3^{}+10^{}=13^{}, \ 3^{7}+10^{}=13^{3}; \nonumber\\
\label{ex-excep-abc} 2^{5}+3^{}=35^{}, \ 2^{3}+3^{3}=35^{};\\ %<== 3+5=2^3, 3^3+5=2^5
2^{}+89^{}=91^{}, \ 2^{13}+89^{}=91^{2};\nonumber\\ 
2^{7}+5^{}=133^{}, \ 2^{3}+5^{3}=133^{}; \nonumber\\ %<== 3+5=2^3, 3+5^3=2^7
2^{8}+3^{}=259^{}, \ 2^{4}+3^{5}=259^{}; \nonumber\\ %<== 3+13=2^4, 3^5+13=2^8
3^{7}+13^{}=2200^{}, \ 3^{}+13^{3}=2200^{}; \nonumber\\ %<== 3+10=13, 3^7+10=13^3
2^{13}+91^{}=8283^{}, \ 2^{}+91^{2}=8283^{}; \nonumber\\ %<== 2+89=91, 2^{13}+89=91^2
2^{}+(2^r-1)^{}={2^r+1}^{}, \ 2^{r+2}+(2^r-1)^{2}=(2^r+1)^{2}.\nonumber
\end{gather}

After the mentioned work of Miyazaki and Pink, they have extensively studied Conjecture \ref{atmost1}, and obtained various finiteness results on it, especially, for cases where the value of $c$ is fixed (cf.~\cite{MiyPi2,MiyPi3}).
In the present article, we complete one of the main purposes in the series of their papers \cite{MiyPi,MiyPi2,MiyPi3}, by establishing a general and satisfactory result on Conjecture \ref{atmost1}, which is stated as follows:

\begin{thm}\label{th-anyc}
For any fixed $c,$ there is at most one solution to equation \eqref{abc}$,$ except for only finitely many pairs of $a$ and $b.$
%Let $c$ be any fixed positive integer greater than $1.$
%Then $N(a,b,c) \le 1,$ except for only finitely many pairs of $a$ and $b.$
\end{thm}

A novelty of this theorem is that it is regarded as a 3-variable generalization of the result of Miyazaki and Pink \cite[Theorem 1]{MiyPi3}, a general and satisfactory contribution to the mentioned conjecture of Bennett on Pillai's type equation, which asserts that for any fixed positive integer $a$ there are only finitely many pairs $(b,c)$ of coprime positive integers with $b>1$ such that the Pillai's type equation $a^x-b^y=c$ has more than one solution in positive integers $x$ and $y$.
It should be remarked that Theorem \ref{th-anyc} is ineffective in the sense that our proof for each $c$ does not provide a way to effectively find exceptional pairs of $a$ and $b$ stated in it.

In what follows, we let $N=N(a,b,c)$ denote the number of solutions $(x,y,z)$ to equation \eqref{abc}.
In the proof of Theorem \ref{th-anyc}, we deal with the following system of two equations:
\begin{eqnarray}\label{abc-sys}
\begin{cases}
\,a^x+b^y=c^z,\\
\,a^X+b^Y=c^Z
\end{cases}
\end{eqnarray}
in positive integers $x,y,z,X,Y$ and $Z$ with $(x,y,z) \ne (X,Y,Z)$.
It is clear that $N \le 1$ if and only if there is no solution to system \eqref{abc-sys}. 
The proof of Theorem \ref{th-anyc} is three-fold. %can be divided into the following three parts. 
The first is to give strong restrictions of some of the exponential unknowns of the system.
Indeed, most of those are already given in \cite{MiyPi3}, among others, one of them ensures that all $x,y,X,Y$ should be less than some constant which depends only on $c$ and is effectively computable.
The second is to give further restrictions by using the theorem of Ridout \cite{Rid} which is a generalization of Thue-Siegel-Roth theorem. 
Here we mention that our application of it can replace some of the arguments in \cite[Sec.\,10]{MiyPi3} relying on the abc-conjecture by unconditional ones.
As a result, the problem is reduced to studying the following system of concreate exponential equations:
\begin{eqnarray}\label{abc-sys-x2y1X1Y2}
\begin{cases}
\,a^2+b=c^z,\\
\,a+b^2=c^Z.
\end{cases}
\end{eqnarray}
%Finally, we handle this system. 
Interestingly, there is just an elementary and `complete' treatment of the above system (see Theorem \ref{th-x2y1X1Y2} below), and it leads us to finish the proof.

The organization of this paper is as follows. 
In the next section, we quote several results of Miyazaki and Pink \cite{MiyPi3} with brief explanations on those proofs.
Section \ref{sec-3} is devoted to prove the main part of the proof of Theorem \ref{th-anyc}, where we are conscious of which parts are effective procedures. 
To finish the proof of Theorem \ref{th-anyc}, we resolve system \eqref{abc-sys-x2y1X1Y2} completely in Section \ref{sec-4}.
The remaining sections are devoted to thoroughly discuss how effectively exceptional pairs of $a$ and $b$ stated in Theorem \ref{th-anyc} for each $c$ can be determined.
After seeing a state-of-the-art result in such direction in Section \ref{sec-5}, based on the methods of proving Theorem \ref{th-anyc}, in Section \ref{sec-6} we observe how useful the theory of rational approximation of algebraic irrationals is, which results problems of establishing effective irrationality measures of many rational powers of $c$.
In Section \ref{sec-7}, we attempt to reduce the number of considerations of such irrationals by giving explicit severe restrictions on exponential unknowns in the system under consideration, in particular, relatively sharp upper bounds for them.
In the final section, we propose a problem towards establishing effective version of Theorem \ref{th-anyc} for certain, presumably infinitely many, values of $c$.

%%%%%%%%%%%%%%%%%%%%%%%%%%%
\section{Works of Miyazaki and Pink}\label{sec-2}
%%%%%%%%%%%%%%%%%%%%%%%%%%%

Here and henceforth, we frequently use the Vinogradov notation 
\[
f \ll_{\kappa_1,\kappa_2,\ldots,\kappa_n} g
\]
which means that $|f/g|$ is less than some positive constant depending only on $\kappa_1,\kappa_2,\ldots$ and $\kappa_n$, where we simply write $f \ll g$ if the implied constant is absolutely finite.
Note that the implied constant of each Vinogradov notation appearing in this paper is effectively computable.

Below, we quote three results, from the work of Miyazaki and Pink, which are direct consequences from the contents of \cite[Sec.\,4]{MiyPi3} (cf.~\cite[Lemmas 4.8, 4.10 and 4.12]{MiyPi3}).
Our proof of Theorem \ref{th-anyc} is based on the first of them, and the other two also play important roles in it.
%Below, we quote each of them with a brief explanation on its proof.

\begin{lem}\label{lem-xyXY-finite}
Let $(x,y,z,X,Y,Z)$ be any solution to system \eqref{abc-sys}$.$
Then
\[
\max\{x,y,X,Y\} \ll_{c}\,1.
\]
\end{lem}

The proof of this lemma proceeds as follows.
It suffices to consider the case where $z \le Z$.
Let $p$ be any prime factor of $c$.
By simple arguments with Baker's method in its $p$-adic form, we use \eqref{abc-sys} to find that both $a,b$ have to be close to 1 in $p$-adic sense, more precisely, both the multiplicative orders of $a,b$ modulo some very large divisor of $p^z$ are small. 
This fact is efficient for finding a sharp upper bound for the $p$-adic valuation of $a^X+b^Y\,(=c^Z)$ being at least $Z$.
These bounds together imply the assertion.

\begin{lem}\label{lem-y1X1}
There exists some constant $K_1$ which depends only on $c$ and is effectively computable such that if $\max\{a,b\}>K_1,$ then
\[
\min\{x,y\}=1, \ \ \min\{X,Y\}=1
\]
for any solution $(x,y,z,X,Y,Z)$ to system \eqref{abc-sys}$.$
\end{lem}

For showing the first inequality of this lemma, we rewrite the first equation in \eqref{abc-sys} as
\[
\mathcal X^m+\mathcal Y^n=c^z,
\]
where $\mathcal X=a,\,\mathcal Y=b,\,m=x,\,n=y$.
Note that for given $c$ Lemma \ref{lem-xyXY-finite} says that both $x,y$ can take only finitely many values.
It is known that for each pair $(m,n)$ of positive integers with $m \ge 2, \,n \ge 2$ and $\max\{m,n\} \ge 3$, the greatest prime factor of $\mathcal X^m+\mathcal Y^n$ with nonzero coprime integers $\mathcal X$ and $\mathcal Y$ tends to infinity as $\max\{|\mathcal X|,|\mathcal Y|\} \to \infty$ in effective sense (cf.~\cite[Corollary 8.1]{ShTi}, \cite{Bu}).
We apply this fact for each possible pair $(x,y)$.
As a result, together with another treatment of case $(x,y)=(2,2)$, we obtain $\max\{a,b\} \ll_c 1$ unless $\min\{x,y\}=1$.
Similarly the second asserted inequality is proved using the second equation in \eqref{abc-sys}.

\begin{lem}\label{lem-xneXyneY}
There exists some constant $K_2$ depending only on $c$ such that if $\max\{a,b\}>K_2,$ then
\[
x \ne X, \ \ y \ne Y
\]
for any solution $(x,y,z,X,Y,Z)$ to system \eqref{abc-sys}$.$
\end{lem}

Suppose that $x=X$.
By \eqref{abc-sys}, we have
\[
\mathcal Z^m-\mathcal Z^n=q^{y_1}-q^{y_2},
\]
where $(m,n,q)=(Y,y,c)$ and $(\mathcal Z,y_1,y_2)=(b,Z,z)$.
Bugeaud and Luca \cite[Theorem 4.2]{BuLu} proved as a corollary of their work that for any positive integers $m,n$ and $q$ with $m>n$ and $q>1$ there are only finitely many solutions to the above equation in positive integers $\mathcal Z,y_1$ and $y_2$ with $\mathcal Z>1$ and $\gcd(\mathcal Z,q)=1$.
Taking the fact that $y,Y \ll_c 1$ by Lemma \ref{lem-xyXY-finite} into consideration, we apply their result for each possible pair $(y,Y)$, and it turns out that both $a,b$ have to be less than some constant depending only on $c$.
Case where $y=Y$ is handled similarly.
We should remark that the mentioned result of Bugeaud and Luca is ineffective in the sense that their proof heavily depends on not only a theorem of Ridout (cf.~Proposition \ref{prop3} below) but also Schmidt Subspace Theorem (cf.~\cite{Sc,Sc2}), thereby it provides no way to effectively find $K_2$.
We also emphasize that Subspace Theorem in that proof is used only when $m \ge 3,\,n=1$ and $q$ is composite (cf.~proof of \cite[Proposition 4.1]{MiyPi3}). 

%%%%%%%%%%%%%%%%%%%%%%%%%%%
\section{Proof of Theorem \ref{th-anyc} for main case}\label{sec-3}
%%%%%%%%%%%%%%%%%%%%%%%%%%%

For our purpose of proving Theorem \ref{th-anyc}, we may assume, with the notation in Lemma \ref{lem-y1X1}, that 
\[
\max\{a,b\}>K_1.
\]
Now, suppose that system \eqref{abc-sys} has a solution $(x,y,z,X,Y,Z)$.
We will observe that this leads to the existence of some upper bound for both $a$ and $b$ which depends only on $c$.

Lemmas \ref{lem-xyXY-finite} and \ref{lem-y1X1} say that
\begin{gather*}
x,y,X,Y \ll_{c}\,1\,;\\
\min\{x,y\}=1, \ \ \min\{X,Y\}=1,
\end{gather*}
respectively.
Further, by the symmetry of $a$ and $b$, without loss of generality, we may assume 
\[
y=1.
\]
To sum up, %we have
\begin{eqnarray}
&a^x+b=c^z,\label{1st}\\
&a^X+b^Y=c^Z \label{2nd}
\end{eqnarray}
with
\[
x \ll_c\,1, \ \ \min\{X,Y\}=1, \ \ \max\{X,Y\} \ll_c\,1.
\]

\begin{lem}\label{Lem-1}
Assume that $X=1.$
Then there exists some constant $K_3$ which depends only on $c$ and is effectively computable such that if $\max\{a,b\}>K_3,$ then the following hold.%\,{\rm:}
\begin{itemize}
\item[{\rm (i)}] $a^x>b.$
\item[{\rm (ii)}] $a<b^Y.$
\end{itemize}
\end{lem}

\begin{proof}
There is no loss of generality in assuming $z \le Z$.
Clearly, $Y>1$.

Reducing \eqref{1st}, \eqref{2nd} modulo $c^z$ gives
\[
a^x \equiv -b \mod{c^z}, \quad a \equiv -b^Y \mod{c^z},
\]
respectively.
Then
\[
a^{x Y} \equiv (-b)^Y \equiv (-1)^Y b^Y \equiv (-1)^{Y+1} a\mod{c^z}.
\]
Since $\gcd(a,c)=1$, one has
\[
a^{x Y-1} \equiv (-1)^{Y+1} \mod{c^z}.
\]
Similarly, $b^{x Y-1} \equiv (-1)^{x+1} \pmod{c^z}$.
Since $xY \ge Y>1$, these two congruences imply that 
\begin{equation}\label{congcz-ineq}
a^{x Y-1} \ge c^z-1, \quad b^{x Y-1} \ge c^z-1,
\end{equation}
respectively.
Assertion (ii) easily follows from the above second inequality as $a^x=c^z-b<c^z-1 \le b^{x Y-1}< b^{xY}$.
Therefore, it remains to show (i).

By equation \eqref{1st}, it suffices to show that $c^z>2b$.
We shall derive a nontrivial lower bound for $c^z/b$.
Since $b^Y<c^Z$ by equation \eqref{2nd}, it follows that
\[
\Big(\frac{c^z}{b}\Big)^{\!Y}=\frac{c^{Y z}}{b^Y}>\frac{c^{Y z}}{c^Z}=c^{Y z - Z},
\]
so that
\begin{equation} \label{ineq-1}
\log \Bigr(\frac{c^z}{b}\Bigr) > \frac{\log c}{Y} \cdot (Y z - Z).
\end{equation}
On the other hand, from equations \eqref{1st} and \eqref{2nd}, observe that 
\[
c^{Y z}=(a^x+b)^Y>a^{x Y} + b^Y > a+b^Y=c^Z.
\]
Thus 
\[
Y z>Z.
\] 
Further, reducing equations \eqref{1st}, \eqref{2nd} modulo $a$ gives
\[
b \equiv c^z \mod{a}, \quad b^Y \equiv c^Z \mod{a},
\]
respectively.
Then 
\[
c^{Y z} \equiv b^Y \equiv c^Z \mod{a},
\] 
so that
\[
c^{Yz -Z} \equiv 1 \mod{a}.
\]
Since $Yz-Z>0$, one has $c^{Yz -Z} \ge a+1$, in particular, 
\begin{equation} \label{ineq-2}
Yz -Z> \frac{\log a}{\log c}.
\end{equation}
One combines \eqref{ineq-1} with \eqref{ineq-2} to find that 
\[
\log \Bigr(\frac{c^z}{b}\Bigr) > \frac{\log c}{Y} \cdot \frac{\log a}{\log c}=\log a^{1/Y},
\]
namely, 
\[
\frac{c^z}{b} > a^{1/Y}.
\]

Finally, suppose that $c^z \le 2b$.
Then $a^{1/Y}<2$, whence $a<2^Y \ll_{c}1$ as $Y \ll_{c}1$. 
Further, since $c^z \le a^{x Y -1}+1$ by the first inequality in \eqref{congcz-ineq}, and $x Y \ll_{c}1$, one obtains $c^z \ll_{c}1$, leading to $b<c^z \ll_{c}1$.
This completes the proof.
\end{proof}

In what follows, we assume, with the notation in Lemmas \ref{lem-xneXyneY} and \ref{Lem-1}, that 
\[
\max\{a,b\}>K,
\]
where $K=\max\{K_1,K_2,K_3\}$, which is some constant depending only on $c$.
Then Lemma \ref{lem-xneXyneY} says that $x \ne X$ and $Y \ne 1$, whence
\[
x \ge 2, \ \ X=1, \ \ Y \ge 2.
\]
In particular, equation \eqref{2nd} becomes
\begin{equation}\label{2nd-2}
a+b^Y=c^Z. 
\end{equation}

\begin{lem}\label{Lem-2}
For any positive number $\varepsilon,$ there exist some positive constants $\mathcal C_1$ and $\mathcal C_2$ each of which depends only on $c$ and $\varepsilon$ such that each of the following inequalities holds.%\,{\rm:}
\begin{itemize}
\item[{\rm (i)}] $b>\mathcal C_1 \,a^{x-1-\varepsilon}.$
\item[{\rm (ii)}] $a>\mathcal C_2 \,b^{Y-1-\varepsilon}.$
\end{itemize}
\end{lem}

To show this lemma, we rely on the following which is a particular case of a generalization of Thue-Siegel-Roth theorem due to Ridout \cite{Rid}.

\begin{prop}\label{prop3}
Let $\xi$ be a nonzero algebraic number.
Let $\mathcal S$ be a non-empty set of finitely many primes.
Then, for every positive number $\varepsilon,$ there exists some positive constant $C$ depending only on $\xi,\mathcal S$ and $\varepsilon$ such that the inequality
\[
\left|\xi - \frac{p}{q}\right| > \frac{C}{q^{1+\varepsilon}}
\]
holds for all integers $p$ and $q$ with $q>0$ and $p/q \ne \xi$ such that every prime factor of $q$ belongs to $\mathcal S.$
\end{prop}

\begin{proof}[Proof of Lemma $\ref{Lem-2}$]
(i) On equation \eqref{1st}, write
\[
z=Q\,x+R
\]
for some nonnegative integers $Q,R$ with $0 \le R<x$. 
Observe that
\begin{align*}
&b=c^z-a^x=(c^{\,z/x}-a) \cdot \frac{c^z-a^x}{\,c^{\,z/x}-a\,},\\
&c^{\,z/x}-a= c^{\,Q+R/x}-a = c^Q\,\Big(c^{\,R/x}-\frac{a}{c^Q}\Bigr),\\
&\frac{c^z-a^x}{\,c^{\,z/x}-a\,}=\big(\,(c^{\,z/x})^{x-1}+(c^{\,z/x})^{x-2}\,a+\cdots+a^{x-1}\,\big) > a^{x-1}.
\end{align*}
These together imply that %$b \ge a^{x-1}$ if $R=0$, and also that
\begin{equation} \label{ineq-approx}
\frac{b}{c^Q \cdot a^{x-1}} > \xi(c,x,R)-\frac{a}{c^Q} \ \, (>0),
\end{equation}
%for $R>0$, 
where 
\[
\xi(c,x,R)=c^{\,R/x}.
\]
Here we note that the inequality $x \ll_{c}1$ ensures that $x,R$ can take only finitely many values for each $c$.

Now, to find a lower bound for the right-hand side of \eqref{ineq-approx}, we shall apply Proposition \ref{prop3}. 
For each $c$, and for each possible pair $(x,R)$, we set $\xi:=\xi(c,x,R)$ and choose $\mathcal S$ as the set of all primes dividing $c$.
Let $\varepsilon$ be any positive number.
Further, we set $p:=a$ and $q:=c^Q$.
It turns out that there is a constant $C=C(c,\varepsilon)>0$ depending only on $c$ and $\varepsilon$ such that the inequality
\begin{equation} \label{ineq-afterridout}
\biggr| \xi(c,x,R)-\frac{a}{c^Q}\biggr| > \frac{C(c,\varepsilon)}{(c^Q)^{\,1+\varepsilon}}
\end{equation}
holds for all possible pairs $(x,R)$.
Inequalities \eqref{ineq-approx}, \eqref{ineq-afterridout} together imply that
\[
b>C(c,\varepsilon)\,\frac{a^{x-1}}{(c^Q)^{\varepsilon}}.
\]
Observe that $c^Q=c^{\,z/x-R/x} \le c^{z/x}$, and that $c^{z/x}<(2a^x)^{1/x} \le 2a$ as $c^z=a^x+b<2a^x$ by Lemma \ref{Lem-1}\,(i).
Therefore,
\[
b>C(c,\varepsilon)\,\frac{a^{x-1}}{(2a)^\varepsilon}
=\mathcal C_1 \,a^{x-1-\varepsilon}%,
\]
with $\mathcal C_1=C(c,\varepsilon)/2^{\varepsilon}$.
\par
(ii) The assertion is proved similarly by observing equation \eqref{2nd-2} with the inequality $Y \ll_{c}1$ and the inequality $c^Z<2b^Y$ following from Lemma \ref{Lem-1}\,(ii). 
\end{proof}

\begin{lem}\label{Lem-3}
If $x>2$ or $Y>2,$ then $\max\{a,b\}$ is less than some constant depending only on $c.$
\end{lem}

\begin{proof}
We use Lemma \ref{Lem-2} for any fixed positive number $\varepsilon$ with $\varepsilon<1$.
Since $x-1-\varepsilon \ge 1-\varepsilon>0$ as $x \ge 2$, one finds that
\begin{equation}\label{ineq-b-fin}
b>\mathcal C_1\,a^{x-1-\varepsilon}
>\mathcal C_1\,(\,\mathcal C_2\,b^{Y-1-\varepsilon}\,)^{x-1-\varepsilon}
=\mathcal C_1\,{\mathcal C_2}^{x-1-\varepsilon}\,b^{\,(x-1-\varepsilon)(Y-1-\varepsilon)}.
\end{equation}
Now, we set $\varepsilon:=1/3$. 
Then
\[
\min\{x-1-\varepsilon,\,Y-1-\varepsilon\} \ge 2-1-1/3=2/3.
\]
If $x>2$ or $Y>2$, then
\[
(x-1-\varepsilon)(Y-1-\varepsilon) \ge 2/3 \cdot (1+2/3)=10/9,
\]
thereby inequality \eqref{ineq-b-fin} gives
\[
b>\mathcal C_3 \,b^{10/9},
\]
where $\mathcal C_3 \ (=\mathcal C_1\,{\mathcal C_2}^{2/3})$ is a positive constant depending only on $c$. 
Therefore, $b$ is less than some constant depending only on $c$. 
Also the assertion on $a$ follows as $a^{2/3} \le a^{x-1-\varepsilon}<b/\mathcal C_1$.
\end{proof}

To sum up, we have established the following:

\begin{lem}\label{summary-xgt2orYgt2}
If $\max\{a,b\}$ exceeds some constant depending only on $c,$ then $(x,y,X,Y)=(2,1,1,2)$ or $(1,2,2,1)$ for any solution $(x,y,z,X,Y,Z)$ to system \eqref{abc-sys}$.$
\end{lem}

According to this lemma, to complete the proof of Theorem \ref{th-anyc}, it suffices to handle system \eqref{abc-sys-x2y1X1Y2}. 
Actually, it is possible to solve it completely, as follows:

\begin{thm}\label{th-x2y1X1Y2}
All solutions to the system of simultaneous equations
\[
a^2+b=c^z, \ \ a+b^2=c^Z
\] 
in positive integers $a,b,c,z$ and $Z$ such that $\min\{a,b,c\}>1$ and $\gcd(a,b)=1$ are given by $(a,b,c,z,Z)=(2,5,3,2,3)$ or $(5,2,3,3,2).$ 
\end{thm}

%%%%%%%%%%%%%%%%%%%%%%%%%%%
\section{Proof of Theorem \ref{th-x2y1X1Y2}}\label{sec-4}
%%%%%%%%%%%%%%%%%%%%%%%%%%%

%\begin{proof}%[Proof of Theorem $\ref{th-x2y1X1Y2}$]
Let $(a,b,c,z,Z)$ be a solution to the system under consideration.
Then
\begin{eqnarray}
&a^2+b=c^z,\label{1st-x2y1X1Y2}\\
&a+b^2=c^Z. \label{2nd-x2y1X1Y2}
\end{eqnarray}
It is easy to see that $a,b,c$ are pairwise coprime.
There is no loss of generality in assuming $a<b$.
Clearly, $z<Z$. 

By \eqref{1st-x2y1X1Y2} and \eqref{2nd-x2y1X1Y2}, %observe that
\[
c^{2z}=(a^2+b)^2=a^4+2a^2b+b^2>a+b^2=c^Z.
\]
Thus
\begin{equation}\label{ineq-2zZ}
2z > Z.
\end{equation}

Since $c^{2z}=a^4+2a^2b+b^2$, and $b^2 \equiv -a \pmod{c^Z}$ by \eqref{2nd-x2y1X1Y2}, 
it follows from \eqref{ineq-2zZ} that $a^4+2a^2b-a \equiv 0 \pmod{c^Z},$ so that 
\[
a\,(a^3+2ab-1) \equiv 0 \mod{c^Z}.
\]
Since $\gcd(a,c)=1$, one finds that
%\begin{equation}\label{rel-f}
\[
a^3+2ab-1=f c^Z
\]
%\end{equation}
for some $f \in \mathbb N$ coprime to $a$.
Further, since $c^Z=a+b^2$, it follows that
%\begin{align} \label{rel-b-af}
\[
f b^2 - 2a b - a^3+f a+1=0.
\]
%\end{align}
Regarding this as a quadratic equation for $b$ with $b>a$, one has
\begin{align} \label{rel-b-af-2}
b%=\frac{1}{f} \,(\,a+\sqrt{(-a)^2+f(a^3-af-1)}\,)
=\frac{a+\sqrt{D}}{f},
\end{align}
where $D$ is a perfect square given by
\[
D=f a^3+a^2-f^2 a-f.
\]
Observe that $b$ is determined by $a$ and $f$, and that $D \equiv -f \pmod{a}$.

Assume that $f=1$.
By \eqref{rel-b-af-2},
\[
b=a+\sqrt{a^3+a^2-a-1}=a+(a+1)\sqrt{a-1}.
\]
Since $a-1$ is a square number, one can write 
\[
a=r^2+1
\]
for some $r \in \mathbb N$.
Observe that
\begin{align*}
&b=a+(a+1)r=r^3+r^2+2r+1,\\
&c^z=a^2+b=r^4+r^3+3r^2+2r+2,\\
&c^Z=a+b^2=r^6+2r^5+5r^4+6r^3+7r^2+4r+2.%,\\
%&.
\end{align*}
Further, recalling that $z<Z<2z$, one calculates that
\[
c^{Z-z}=r^2+r+1, \ \ c^{2z-Z}=r^2+2.
\]
Euclidean algorithm for these two integers reveals that
\[
\gcd(c^{Z-z},c^{2z-Z}) \,\mid \, 3.
\]
This implies that $\min\{c^{Z-z},c^{2z-Z}\}=3$.
Thus $c=3,r=1$, so that $a=2,b=5,z=2,Z=3$.
Therefore, in what follows, we suppose that
\begin{equation}\label{ineq-f>_2}
f \ge 2.
\end{equation}
We will observe that this leads to a contradiction.

Since $a^2 \equiv -b \pmod{c^z}$ and $a \equiv -b^2 \pmod{c^z}$, one has 
\[
ab \equiv 1 \mod{c^z}.
\]
Thus
\[
ab-1=g c^z=g(a^2+b)
\]
for some $g \in \mathbb N$.
Reducing this modulo $a$ gives
\[
-1 \equiv g b \mod{a}.
\]
On the other hand, reducing \eqref{rel-b-af-2} modulo $a$ implies
\[
fb \equiv \sqrt{D} \mod{a}.
\]
Since $\gcd(b,a)=1$, one combines these two congruences modulo $a$ to easily find that 
\[
-f \equiv g \sqrt{D} \mod{a}.
\]
Squaring both sides of the above congruence gives $f^2 \equiv g^2 D \equiv g^2 \cdot (-f) \pmod{a}$, so that $f \equiv -g^2 \pmod{a}$ as $\gcd(f,a)=1$.
In particular, 
\begin{equation}\label{ineq-fgg>_a}
f+g^2 \ge a.
\end{equation}

Since 
\[
f =\frac{a^3+2ab-1}{c^Z}=\frac{a(a^2+b)+(ab-1)}{c^Z}=\frac{ac^z+gc^z}{c^Z}=\frac{a+g}{c^{Z-z}},
\]
one uses \eqref{ineq-fgg>_a} to find that
\[
\frac{a+g}{c^{Z-z}}+g^2 \ge a,
\]
whence
\[
\biggr(1-\frac{1}{c^{Z-z}}\biggr)\,a \le \frac{g}{c^{Z-z}}+g^2.
\]
Further, since
\[
g=\frac{ab-1}{a^2+b}
<\frac{ab}{a^2}=\frac{b}{a}
=\frac{a+\sqrt{D}}{fa}
<\frac{a+\sqrt{f a^3+a^2}}{fa}
=\frac{1+\sqrt{fa+1}}{f}
\]
by \eqref{rel-b-af-2}, it follows that
\begin{equation}\label{ineq-a}
\biggr(1-\frac{1}{c^{Z-z}}\biggr)\,a < \frac{g_u}{c^{Z-z}}+{g_u}^2,
\end{equation}
where 
\[
%g_u={g_u}(a,f):=\frac{1}{f}+\sqrt{\frac{a}{f}+\frac{1}{f^2}-\frac{1}{a}-\frac{1}{fa^2}}.
g_u={g_u}(a,f):=\frac{1}{f}+\sqrt{\frac{a}{f}+\frac{1}{f^2}}.
\]
Inequality \eqref{ineq-a} implies that $a$ is absolutely finite, whenever
\[
1-\frac{1}{c^{Z-z}}>\frac{1}{f}.
\]

Indeed, suppose that $c^{Z-z}>2$.
By \eqref{ineq-a},
\[
\frac{2}{3}\,a \le \frac{g_u}{3}+{g_u}^2.
\]
Since $f \ge 2$ by \eqref{ineq-f>_2}, the above displayed inequality, together with the fact that $D$ is a perfect square and $b$ is an integer given by formula \eqref{rel-b-af-2}, implies that $a=3$ with $f=2$, and $b=5$.
However, equation \eqref{2nd-x2y1X1Y2} does not hold as $Z>1$.

Finally, suppose that $c^{Z-z}=2$.
Then $c=2$ and $Z=z+1$.
From equations \eqref{1st-x2y1X1Y2} and \eqref{2nd-x2y1X1Y2},
\begin{align*}
a+b^2 - (a^2+b) &= (b-a)(b+a-1) \\
&=2^{z+1}-2^z= 2^z.
\end{align*}
Since $a,b$ are odd, the factor $b+a-1$ equals $1$, that is, $b+a=2$, which clearly contradicts the premise that $\min\{a,b\}>1$.
This completes the proof of Theorem \ref{th-x2y1X1Y2}.
%\end{proof}

To sum up, Lemma \ref{summary-xgt2orYgt2} and Theorem \ref{th-x2y1X1Y2} together complete the proof of Theorem \ref{th-anyc}.

In the remaining sections, we consider effective version of Theorem \ref{th-anyc}, namely, we discuss how effectively exceptional pairs of $a$ and $b$ stated in the theorem for each $c$ can be determined.

%%%%%%%%%%%%%%%%%%%%%%%%%%%
\section{Effective version of Theorem \ref{th-anyc} for certain cases}\label{sec-5}
%%%%%%%%%%%%%%%%%%%%%%%%%%%

Here, we just introduce a state-of-the-art result on effective version of Theorem \ref{th-anyc}.
To the best of our knowledge, we have the following:

\begin{prop}\label{prop-stateoftheart}
Let $c$ be any fixed positive integer satisfying at least one of the following conditions\,{\rm :}
\begin{itemize}
\item[\rm (i)] $\max\{\,2^{\nu_2(c)},3^{\nu_3(c)}\,\}>\sqrt{c},$ where $\nu_{p}$ denotes the $p$-adic valuation\,$;$
\item[\rm (ii)] $c$ is a prime of the form $2^r+1$ with some positive integer $r$\,$;$
\item[\rm (iii)] $c=13.$
\end{itemize}
Then $N(a,b,c) \le 1,$ except for only finitely many pairs of $a$ and $b,$ all of which are effectively determined.
\end{prop}

Cases (i), (ii) and (iii) in this proposition are handled by \cite[Corollary 2]{MiyPi2}, the proof of \cite[Theorem 3]{MiyPi2} and \cite[Theorem 3]{MiyPi3}, respectively.
The values of $c$ treated in Proposition \ref{prop-stateoftheart}, not being perfect powers, are in ascending order as follows:
\[
c=2, 3, 5, 6, 12, 13, 17, 18, 24, 40, 45, 48, 54, 56, 63,\ldots,257,\ldots,65537,\ldots
\]
We note that among the above cases Conjecture \ref{atmost1} is completely solved for each of the following cases: 
\begin{itemize}
\item[$\bullet$] $c=2$ by Scott \cite[Theorem 6; $p=2$]{Sc_JNT_93}\,$;$
\item[$\bullet$] $c=6$ or $c \in \{3,5,17,257,65537\}$ (Fermat primes found so far) by Miyazaki and Pink \cite{MiyPi2}\,$;$
\item[$\bullet$] $c=13$ by Miyazaki and Pink \cite{MiyPi3}.
\end{itemize}
Further, related to case (i), we mention that it is possible to prove Conjecture \ref{atmost1} if $\max\{2^{\nu_2(c)},3^{\nu_3(c)}\}>c^{1-\epsilon}$ for a small constant $\epsilon>0$ (cf.~proof of \cite[Theorem 2]{MiyPi2}).

%%%%%%%%%%%%%%%%%%%%%%%%%%%
\section{Towards effective version of Theorem \ref{th-anyc} : rational approximations of algebraic irrationals}\label{sec-6}
%%%%%%%%%%%%%%%%%%%%%%%%%%%

In this section we shall discuss how effective Theorem \ref{th-anyc} can be.
Namely, we attempt to find a constant $K$ which depends only on $c$ and is effectively computable such that $N(a,b,c) \le 1$ whenever $\max\{a,b\}>K$.
For this, we examine the parts in our proof of Theorem \ref{th-anyc} which relies on rational approximations of algebraic irrationals.
One of the conclusions of this section is that the effective version of Theorem \ref{th-anyc} for prime values of $c$ is obtained if all Problems \ref{prob1}, \ref{prob2} and \ref{prob3} below are solved.

For our purpose, as discussed in Section \ref{sec-3}, we may assume that we have
\begin{eqnarray}
&a^x+b=c^z,\label{eff-1st-y1}\\
&a^X+b^Y=c^Z \label{eff-2nd-y1}
\end{eqnarray}
with
\[
x \ll_c\,1, \ \ \min\{X,Y\}=1, \ \ 2 \le \max\{X,Y\} \ll_c\,1, \ \ z \le Z.
\]
Below, we always consider the system of equations \eqref{eff-1st-y1} and \eqref{eff-2nd-y1}.
There are two cases according to whether $X=1$ or $Y=1$.
We remark that one can omit considering the case where $a$ or $b$ is congruent to 1 or $-1$ modulo $c$ by \cite[Corollary 1]{MiyPi2}.

\subsection{Case where $X=1$}
We distinguish two cases according to whether $x \ge 2$ or $x=1$.

\subsubsection{Case where $X=1$ and $x \ge 2$}

Recall that $x,Y \ll_c 1$.
We are in the same situation as that of Lemma \ref{Lem-2}. 
For each possible pair $(x,Y)$, and for each of the pairs $(r,R)$ of nonnegative integers with $r<x$ and $R<Y$, we put
\[
\xi_1={\xi_1}(c,x,r)=c^{\,r/x}, \quad
\xi_2={\xi_2}(c,Y,R)=c^{\,R/Y}. 
\]
Under this setting, the problem corresponding to this case is the following:

\begin{prob}\label{prob1}
\rm
For each quadruple $(x,Y,r,R)$, find positive constants $\epsilon_1$ and $\epsilon_2$ with $(x-1-\epsilon_1)(Y-1-\epsilon_2)>1$ and positive constants $C_1=C_1(\xi_1,\epsilon_1)$ and $C_2=C_2(\xi_2,\epsilon_2)$ such that each of the inequalities 
\[
\left| \xi_1 - \frac{p}{q} \right| > \frac{C_1}{q^{x-1-\epsilon_1}} 
\text{ \ and \ }
\left| \xi_2 - \frac{p}{q} \right| > \frac{C_2}{q^{Y-1-\epsilon_2}}
\]
holds for any positive integers $p$ and $q$ with the property that $q$ is a power of $c$.
\end{prob}

\begin{rem}\label{rem-coprime}
\rm
By the well-known theorem of Liouville \cite{Li}, for Problem \ref{prob1}, it is essential to consider the case where $\gcd(r,x)=1$ and $\gcd(R,Y)=1$.
A similar remark is given to each of Problems \ref{prob2} and \ref{prob3} below.
\end{rem}

\subsubsection{Case where $X=1$ and $x=1$}\label{sec-X1x1}
As referred to in the explanation on the proof of Lemma \ref{lem-xneXyneY}, this case is handled similarly to the previous one, whenever $c$ is a prime.
More precisely, according to the proof of \cite[Proposition 4.1]{MiyPi3}, with a somewhat nontrivial consideration to show $Y>2$, it turns out that the problem corresponding to this case is the following:
\begin{prob}\label{prob2}
\rm
Assume that $c$ is a prime.
For each pair $(Y,R)$ with $Y \ge 3$, find a positive constant $\epsilon_2$ with $\epsilon_2<\frac{2Y-4}{Y^2-3Y+3}$ and a positive constant $C_2=C_2(\xi_2,\epsilon_2)$ such that the inequality
\[ 
\left| \xi_2 - \frac{p}{q} \right| > \frac{C_2}{q^{1+\epsilon_2}}
\]
holds for any positive integers $p$ and $q$ with the property that $q$ is a power of $c$.
\end{prob}

\subsection{Case where $Y=1$}
The case where $x=1$ can be handled in the same way to that of Section \ref{sec-X1x1}.
For case $x \ge 2$, we can follow the proof of \cite[Proposition 4.1]{MiyPi3}, with a somewhat nontrivial consideration to show $X>2x$, and it turns out that the problem corresponding to this case is the following:
\begin{prob}\label{prob3}
\rm
For each possible $x$ with $x \ge 2$ and for each nonnegative integer $r$ with $r<x$, find a positive constant $\epsilon$ and a positive constant $C=C(c,x,r,\epsilon)$ such that the inequality
\[ 
\left| c^{\,r/x} - \frac{p}{q} \right| > \frac{C}{q^{x-\epsilon}}
\]
holds for any positive integers $p$ and $q$ with the property that $q$ is a power of $c$.
\end{prob}

\begin{rem}\label{rem-feldman}
\rm
By the result of Fel'dman \cite{Fe}, Problem \ref{prob3} can be solved if the algebraic number $c^{\,r/x}$ is of degree at least 3 (see also \cite{BiBu}).
\end{rem}

%%%%%%%%%%%%%%%%%%%%%%%%%%%
\section{Towards effective version of Theorem \ref{th-anyc} : restricting possible values of unknowns}\label{sec-7}
%%%%%%%%%%%%%%%%%%%%%%%%%%%

Continuously, in this section, we always consider the system of equations \eqref{eff-1st-y1} and \eqref{eff-2nd-y1}.
By the contents of the previous section, to make the proof of Theorem \ref{th-anyc} effective, it is clearly important to restrict possible values of $x$ and $\max\{X,Y\}$. 
Below, we give a number of such restrictions.
In what follows, we put 
\[
\Delta=|xY-X|.
\]
We know that $\Delta>0$ (cf.~\cite[Lemma 3.3]{HuLe}).

\subsection{Restricting possible values of $x$}
The main aim of this subsection is to give reasonable upper bounds for $x$ by elementary methods. 
We begin by preparing two notation. 

\begin{defi}\rm
For any positive integer $M$, and an integer $\mathcal A$ coprime to $M$, we define $e_M(\mathcal A)$ to be the least positive integer $e$ such that $\mathcal A^e$ is congruent to $1$ or $-1$ modulo $M$.
\end{defi}

Note that $e_M(\mathcal A)$ is a divisor of $\varphi(M)$, in particular, $e_M(\mathcal A)<M$.

\begin{defi}\rm
For a finite set $S$ of prime numbers, we let $\mathcal A[S]$ denote the $S$-{\it part} of a nonzero integer $\mathcal A$, namely,
\[
\mathcal A[S]=\prod_{p \in S} \,p^{\,\nu_p(\mathcal A)}.%,
\]
For simplicity and convenience, we write $\mathcal A[\,\{p\}\,]=\mathcal A[\,p\,]$ for any prime $p$, and $\mathcal A[\varnothing]=1$, respectively.
\end{defi}

\begin{prop}\label{prop-x-bounds}
If $\max\{a,b\}$ exceeds some constant which depends only on $c$ and is effectively computable, then
\[
x \le \displaystyle \biggl \lfloor \frac{e_{c'}(a) \log c}{\log c[R]}\biggl \rfloor
\]
for any positive divisor $c'$ of $c$ with $c'>1,$ where
\[
R=R(c')=
\begin{cases}
\, S' & \text{ if $c'=2$ or $c' \not\equiv 2 \pmod{4}$},\\
\,S' \setminus \{2\} & \text{ if $c'>2$ and $c' \equiv 2 \pmod{4}$},
\end{cases}
\]
with $S'$ the set of all prime factors of $c'.$
%\end{itemize}
\end{prop}

To show this proposition, we use the following two lemmas, the former of which is an almost direct consequence from the contents of \cite[Sec.\,4]{MiyPi3} (cf.~\cite[Lemma 4.4]{MiyPi3}), and the latter one is well-known in $p$-adic calculations.

\begin{lem}\label{lem-Delta}
Let $c'$ be a positive divisor of $c$ with $c'>1.$
Put 
\[
E'=\lcm(\,e_{c'}(a),e_{c'}(b)\,).
\]
Then the following hold.
\begin{itemize}
\item[\rm (i)]
$\Delta \equiv 0 \pmod{E'}.$
\item[\rm (ii)]
$\gcd(a^{E'}-\delta_a,\,b^{E'}-\delta_b) \cdot \Delta/E' \equiv 0 \pmod{{c'}^z},$ where the number $\delta_h$ with $\delta_h \in \{1,-1\}$ for each $h \in \{a,b\}$ is defined by the following congruence\,$:$
\[ \begin{cases}
\,h^{\,e_{c'}(h)} \equiv \delta_h \mod{c'} & \text{if $c'>2$},\\
\,h \equiv \delta_h \mod{4} & \text{if $c'=2$}.
\end{cases}\]
\end{itemize}
\end{lem}

\begin{lem}\label{padic-lemma}
Let $p$ be a prime number.
Let $U$ and $V$ be relatively prime nonzero integers.
Assume that
\[\begin{cases}
\, U \equiv V \pmod{p} & \text{if $p \ne 2$},\\
\, U \equiv V \pmod{4} & \text{if $p=2$}.
\end{cases}\]
Then
\[
\nu_p(U^N-V^N)=\nu_p(U-V)+\nu_p(N)
\]
for any positive integer $N.$
\end{lem}

\begin{proof}[Proof of Proposition $\ref{prop-x-bounds}$]
The proof proceeds along similar lines to those of \cite[Lemmas 4.4 and 4.11]{MiyPi3}.
First, consider the case where $c'>2$.
Put $e=e_{c'}(a)$.
Then 
\begin{equation}\label{cong1}
a^e \equiv \delta \mod{c'}
\end{equation}
for some $\delta \in \{1,-1\}$.
On the other hand, similarly to the proof of Lemma \ref{Lem-1},
\begin{equation}\label{cong2}
a^{\Delta} \equiv (-1)^{Y+1} \mod{c^z}.
\end{equation}
In particular, 
\begin{equation}\label{cong3}
a^{\Delta} \equiv (-1)^{Y+1} \mod{c'}.
\end{equation}
Since $\Delta \equiv 0 \pmod{e}$ by Lemma \ref{lem-Delta}\,(i), and congruence \eqref{cong3} on the modulus $c'$ ($>\!2$) is also obtained by raising both sides of congruence \eqref{cong1} to the $\Delta/e$-th power, it turns out that $\delta^{\,\Delta/e}=(-1)^{Y+1}$.
Thus, \eqref{cong2} becomes
\[
a^{\Delta} \equiv \delta^{\,\Delta/e} \mod{c^z}.
\]

Let $p$ be any prime factor of $c'$.
By the last congruence above,
\[
\nu_p(\,a^{\Delta}-\delta^{\,\Delta/e}\,) \ge \nu_p(c) \cdot z.
\]
From \eqref{cong1}, observe that $a^e \equiv \delta \pmod{p}$ and also that $a^e \equiv \delta \pmod{4}$ if $p=2$ and $c' \equiv 0 \pmod{4}$. 
Based on these properties, to calculate the left-hand side of the above displayed inequality, one applies Lemma \ref{padic-lemma} for $(U,V)=(a^e,\delta)$ and $N=\Delta/e$.
It turns out that
\[
\nu_{q}(a^e-\delta)+\nu_{q}(\Delta/e) \ge \nu_{q}(c) \cdot z,
\]
where $q$ is any prime factor of $c'$ if $c' \not\equiv 2 \pmod{4}$, and $q$ is any odd prime factor of $c'$ if $c' \equiv 2 \pmod{4}$.
To sum up, 
\[
(a^e-\delta) \cdot \Delta/e \equiv 0 \mod{(c[R])^z}.
\]
%where $R=R(c,c')$.
In particular, 
\[
(c[R])^z \le (a^e+1) \cdot \Delta/e.
\]
Since $a<c^{z/x},$ and $\Delta=|xY-X| \ll_{c}1$, the above inequality implies that
\[
(c[R])^z < c_0 \cdot (c^e)^{z/x}
\]
for some constant $c_0>0$ which depends only on $c$ and is effectively computable.
Raising both sides of the above inequality to the $x/z$-th power gives
\[
(c[R])^x < ({c_0}^x)^{1/z} \cdot c^e.
\]
Observe that ${c_0}^x \ll_c 1$ as $x \ll_{c}1$ and that $e=e_{c'}(a)<c' \ll_c 1$.
It follows from the above displayed inequality that if $z \ge z_0$, where $z_0=\big\lceil \frac{x\log {c_0}}{\log (1+1/c^e)} \bigr\rceil$ with $z_0 \ll_{c}1$, then
\[
(c[R])^x \le c^e.
\]
This gives the asserted inequality.
While if $z<z_0$, then $\max\{a,b\}<c^z<c^{z_0} \ll_c 1$.
This completes the proof for $c'>2$.
For $c'=2$, since one may assume that $z>1$, one can proceed similarly by replacing the modulus $c'$ of each of congruences \eqref{cong1} and \eqref{cong3} by $4$.
\end{proof}

\begin{ex}\label{ex-xleE}\rm
In Proposition \ref{prop-x-bounds}, taking $c'=c$, we have
\[
x \le 
\begin{cases}
\,e_{c}(a) & \text{ if $c=2$ or $c \not\equiv 2 \pmod{4}$}\,;\\
\,\displaystyle \bigg \lfloor \frac{\log c}{\,\log (c/2)} \cdot e_{c}(a) \bigg \rfloor & \text{ if $c>2$ and $c \equiv 2 \pmod{4}$},
\end{cases}
\]
whenever $\max\{a,b\}$ exceeds some constant which depends only on $c$ and is effectively computable.
\end{ex}

\subsubsection{Application of Proposition $\ref{prop-x-bounds}$}\label{sec-xsieve}

Here we shall present how useful Proposition \ref{prop-x-bounds} is for restricting the value of $x$ for some special case.
For this, the following result works efficiently, which is an almost direct consequence from the proof of \cite[Theorem 2]{MiyPi2}.

\begin{lem}\label{Lem-MiyPi2}
Let $S$ be a $($possibly empty$)$ set of odd prime factors of $c.$
Define $M_S$ and $c_S$ as either
\begin{alignat*}{5}
&M_S=\prod_{p \in S}p, & \ \ \ &c_S=\max\{\,c[S],\,c[2]\,\} & \ \ \ \text{or}\\
&M_S=4\prod_{p \in S}p, & & c_S=c\bigr[S \cup \{2\}\bigr].
\end{alignat*}
Assume that $a \equiv \pm 1 \pmod{M_S}$ and $c_S>\sqrt{c}.$
Then the system of equations \eqref{eff-1st-y1} and \eqref{eff-2nd-y1} has no solution, whenever $\max\{a,b\}$ exceeds some constant which depends only on $c$ and is effectively computable.
\end{lem}

In what follows, we consider the case where $X=1$ (this condition can be replaced by $\gcd(\,X,\varphi(c)\,)=1$).
Also we assume that $\max\{a,b\}$ is so large that both the conclusions of Proposition \ref{prop-x-bounds} and Lemma \ref{Lem-MiyPi2} hold.
Proposition \ref{prop-x-bounds} and Lemma \ref{lem-Delta}\,(i) say that
\begin{equation} \label{bounds-xyXY}
x \le \displaystyle \biggl \lfloor \frac{e_{c'}(a) \log c}{\log c[R]}\biggl \rfloor, \quad
x Y \equiv 1 \mod{e_{c'}(a)},
\end{equation}
respectively, where $c'>1$ is any divisor of $c$.
Note that the above congruence condition implies that $x$ is coprime to $e_{c'}(a)$. 
Below, by using restrictions \eqref{bounds-xyXY}, for given $c$ we sieve the value of $x$ by the following strategy, and find a list composed of all possible values of $x$ (which always includes 1).

\smallskip

Take any divisor $c'$ of $c$ with $c'>1$.
Take any integer $a_0$ such that 
\begin{equation} \label{a0-condition}
1 \le a_0 \le \Big\lfloor \,\frac{c'}{2}\, \Big\rfloor \text{ \, and \,} \gcd(a_0,c')=1.
\end{equation}
This $a_0$ corresponds to the congruence class of either $a$ or $-a$ modulo $c'$.
Write $e=e_{c'}(a_0)$.
Let $I_{c',\,a_0}$ be the set of all $t \in \N$ such that
\[
t \le \displaystyle \biggl \lfloor \frac{e \log c}{\log c[R]}\biggl \rfloor \text{ \, and \,} \gcd(t,e)=1.
\]
From \eqref{bounds-xyXY} observe that $I_{c',\,a_0}$ is a list of all possible values of $x$ corresponding to the choice of $c'$ and the case where $a \equiv \pm\,a_0 \pmod{c'}$.
Let $I_{c'}$ be the union of all $I_{c',\,a_0}$, that is, 
\[
I_{c'} = \bigcup_{a_0\text{ with \eqref{a0-condition}}}I_{c',\,a_0}.
\]
Thus, $I_{c'}$ is a list of all possible values of $x$ corresponding to the choice of $c'$.
Here, by applying Lemma \ref{Lem-MiyPi2} with $S$ the set of all odd primes dividing $c'$, we may set $I_{c',\,a_0}=\varnothing$ if $e=1$ and $c_{S}>\sqrt{c}$, where %$c_{S}$ is given by 
\[c_{S}=
\begin{cases}
\,\max\{\,c[S],\,c[2]\,\}
& \text{if $c' \not\equiv 0 \pmod{4}$},\\
\,c\bigr[S \cup \{2\}\bigr]
& \text{if $c' \equiv 0 \pmod{4}$}.
\end{cases}
\]
Let $I$ be the intersection of all $I_{c'}$, that is, 
\[
I = \bigcap_{c' \mid c,\, c'>1}I_{c'}.
\]
This is a list of all possible values of $x$.

\begin{ex}\label{ex-c210}\rm
We demonstrate the above sieve for $c=210$.
We shall take $c'=3,7$ or 15 in turn.

First, let $c'=3$.
Then
\[
x\le \displaystyle \biggl \lfloor \frac{\log 210}{\log 3}\biggl \rfloor=4,
\]
so that $I_{3}=\{1,2,3,4\}$.

Second, let $c'=7$. Then
\[
x\le \displaystyle \biggl \lfloor \frac{e \log 210}{\log 7}\biggl \rfloor, \ \ \gcd(x,e)=1
\]
with $e=e_{7}(a)$.
Since $e \in \{1,3\}$, it follows that
\[
\begin{cases}
\,x\le \textstyle \big \lfloor \frac{\log 210}{\log 7}\big \rfloor=2 & \text{if $e=1$},\\
\,x\le \textstyle \big \lfloor \frac{3\log 210}{\log 7}\big \rfloor=8, \ \ 3 \nmid x& \text{if $e=3$},
\end{cases}
\]
so that $I_{7}=\{1,2\} \cup \{1,2,4,5,7,8\} =\{1,2,4,5,7,8\}$.

Third, let $c'=15$. Then
\[
x \le \displaystyle \biggl \lfloor \frac{e \log 210}{\log 15}\biggl \rfloor, \ \ \gcd(x,e)=1
\]
with $e=e_{15}(a)$.
By Lemma \ref{Lem-MiyPi2}, it suffices to consider when $e>1$.
Since $e \in \{2,4\}$, it follows that
\[
\begin{cases}
\,x\le \textstyle \big \lfloor \frac{2\log 210}{\log 15}\big \rfloor=3, \ \ 2 \nmid x& \text{if $e=2$},\\
\,x\le \textstyle \big \lfloor \frac{4\log 210}{\log 15}\big \rfloor=7, \ \ 2 \nmid x& \text{if $e=4$},
\end{cases}
\]
so that $I_{15}=\{1,3\} \cup \{1,3,5,7\}=\{1,3,5,7\}$.

\smallskip
To sum up, $I_3 \cap I_7 \cap I_{15} =\{1\}$, whence $I=\{1\}$, and $x=1$.
\end{ex}

\begin{rem}\label{rem-x1}\rm
Similarly to Example \ref{ex-c210}, the sieve leads to $x=1$ for infinitely many composite values of $c$ each of which has a relatively large divisor being a power of 2, 3 or 5. 
Such examples are those of $c$ with the property (i) in Proposition \ref{prop-stateoftheart}, and Theorem \ref{th-anyc} for this case is effectively proved in \cite[Corollary 2]{MiyPi2} where a certain property of only primes 2 and 3 is used.
\end{rem}

\subsubsection{Case where $c$ is prime}

Here, we just state a refined version of Proposition \ref{prop-x-bounds} for prime values of $c$.

\begin{prop}%\label{c-general-summary}
Assume that $c$ is an odd prime.
Put $E=\lcm(\,e_{c}(a),e_{c}(b)\,).$
Then the following hold. 
\begin{itemize}
\item[$\bullet$]
If either $E=1$ or $E$ is even, then $\max\{a,b\} \ll_c 1.$
\item[$\bullet$]
Assume that $E$ is odd with $E \ge 3.$
If $\max\{a,b\}$ exceeds some constant which depends only on $c$ and is effectively computable, then $x \le E-2.$
\end{itemize}
\end{prop}

This proposition is an almost direct consequence of the contents of \cite[Sec.\,5]{MiyPi3} (cf.~\cite[Lemma 5.4]{MiyPi3}).
We emphasize that the primality condition on $c$ is essentially used in its proof. 
Indeed, the parities of the exponential unknowns in the left-hand side of each of two equations in the system \eqref{abc-sys} are generally restricted by a result of Scott \cite[Lemma 6]{Sc_JNT_93}, from which methods for studying ternary Diophantine equations (generalized Fermat equations and generalized Lebesgue-Nagell equations, etc.)\! may be also applied in the study of effective version of Theorem \ref{th-anyc} (cf.~proof of \cite[Theorem 3]{MiyPi2} which proves Proposition \ref{prop-stateoftheart} for (ii)).
Further, in the direction of Remark \ref{rem-x1}, we can find many prime values of $c$ for which the sieve described in Section \ref{sec-xsieve} leads to $x=1$ (cf.~proof of \cite[Theorem 3]{MiyPi3} which proves Proposition \ref{prop-stateoftheart} for (iii)).

\subsection{Finding upper bounds for $x$ and $\max\{X,Y\}$}
The main aim of this subsection is to give relatively sharp upper bounds for $x$ and $\max\{X,Y\}$ by Baker's method. 
In what follows, we put
\[
E=\lcm(\,e_{c}(a),e_{c}(b)\,).
\]
We start with the following lemma.

\begin{lem}\label{lem-heggcz}
$a^E \gg_{c} c^z$ and $b^E \gg_{c} c^z.$
\end{lem}

\begin{proof}
Let $h \in \{a,b\}$.
Since $\Delta \ll_{c} 1$, and $(h^E-\delta_h) \cdot \Delta/E \equiv 0 \pmod{c^z}$ by Lemma \ref{lem-Delta}\,(ii) with $c'=c$, it follows that $
c^z \le (h^E+1) \cdot \Delta/E \ll_{c} h^E$.
\end{proof}

There are two cases according to whether $X=1$ or $Y=1$.

\subsubsection{Case where $X=1$}
Consider the following system of two equations:
\begin{eqnarray}
&a^x+b=c^z,\label{eff-1st-y1X1}\\
&a+b^Y=c^Z \label{eff-2nd-y1X1}
\end{eqnarray}
with $Y \ge 2$ and $z \le Z$.
It is not hard to see that $a<b^Y$ (cf.~proof of Lemma \ref{Lem-1}\,(ii)).
%We may assume that $\max\{a,b\}>K_3$.

\begin{prop}\label{prop-systemX1-bounds}
For any number $\delta$ with $0<\delta<1,$ the following hold.
\begin{itemize}
\item[\rm (i)] 
If $a<b^{\delta Y},$ then either 
\[
Y < \min\biggl\{s_0,\,\frac{2520}{1-\delta}\biggr\} \cdot \log c,
\]
or
\[
s_0 \le s<\frac{25.2}{\,1-\delta\,}\,\big( \log (2s+1) +0.38\,\bigr)^2,
\]
where $s=Y/\log c$ and $s_0=\frac{\exp(9.62)-1}{2}=7531.02...\,.$
\item[\rm (ii)] 
Assume that $a \ge b^{\delta Y}.$ 
Then the following hold.
\begin{itemize}
\item[\rm (ii-1)] 
Let $\epsilon$ be any positive number with $\epsilon<1.$
If $\max\{a,b\}$ exceeds some constant which depends only on $c$ and $\epsilon$ and is effectively computable, then
\[
x Y=1+jE
\]
for some positive integer $j$ with $j<(1+\epsilon)/\delta.$
\item[\rm (ii-2)] 
If $\max\{a,b\}$ exceeds some constant which depends only on $c$ and $\delta$ and is effectively computable, then
$x \le 1/\delta.$
% then $z=Z,$ and
%\[
%x <\frac{2520}{\,1-1/Y\,}\log c, \quad z<\frac{\log 2}{\,\bigl(1 -\frac{1}{\delta x}\bigl)\,\log c\,}.
%\]
\end{itemize}
\end{itemize}
\end{prop}

\begin{rem}\rm
On the numerical constant $2520$ (or $25.2$) appearing in Proposition \ref{prop-systemX1-bounds}\,(i), we note that that is an artifice of our proof and can be somewhat reduced. 
The same remark is given to Proposition \ref{prop-systemY1-bounds}\,(i) below.
\end{rem}

\begin{proof}[Proof of Proposition $\ref{prop-systemX1-bounds}$]
(i) Assume that
\begin{equation}\label{alessbdeltaY}
a<b^{\delta Y}.
\end{equation}
Dividing both sides of equation \eqref{eff-2nd-y1X1} by $b^Y$ gives
\[
\frac{a}{b^Y}+1=\frac{c^Z}{b^Y}.
\]
Since
\[
Z \log c - Y \log b=\log \frac{\,c^Z\,}{b^Y}=\log\Bigr(1+\frac{a}{b^Y}\Bigr)
<\frac{a}{b^Y},%<\frac{1}{b^{(1-\delta)Y}}, 
\]
it follows from \eqref{alessbdeltaY} that
\[
\log \varLambda<\log \frac{a}{\,b^Y\,}<- (1-\delta)(\log b)\,Y,
\]
where 
\[
\ \ \ \varLambda=Z \log c - Y \log b \ \ (>0).
\]
On the other hand, from a simple application of \cite[Corollary 2;\,$(m,C_2)=(10,25.2)$]{La},
\[
\log \varLambda > -25.2 \,\log c\,\log b\, \biggl(\,\max \biggl \{ \log \Bigr( \frac{Z}{\log b}+\frac{Y}{\log c} \Bigr) +0.38, \,10\biggl\} \,\biggl)^2.
\]
These bounds for $\log \varLambda$ together yield
\[
s <\frac{25.2}{1-\delta}\,\Bigr( \max\Bigr\{\! \log (2s+1) +0.38,\,10 \Bigr\} \Bigr)^2.
\]
This implies the assertion by noting that $\log (2s+1) +0.38 \ge 10$ if and only if $s \ge s_0$.
\par
(ii-1) Assume that
\begin{equation}\label{agebdeltaY}
a \ge b^{\delta Y}.
\end{equation}
Since $a^x<c^z$, and $c^z \ll_{c}b^E$ by Lemma \ref{lem-heggcz}, one has
\[
a<c^{z/x} \ll_{c} b^{E/x}.
\]
This together with \eqref{agebdeltaY} implies that
\[
b^{\,\delta Y-E/x} \ll_{c} 1.
\]
If $\delta\,Y-E/x \ge \epsilon/2$, then $b^{\epsilon/2} \ll_{c} 1$, so that $b \ll_{c,\epsilon} 1$, leading to $\max\{a,b\}<c^z \ll_c b^E \ll_{c,\epsilon} 1$.
Thus we may assume that $\delta\,Y-E/x < \epsilon/2,$ i.e.,
\[
x Y<\frac{1}{\delta}\,\bigg(\frac{\epsilon}{2}\,x+E\bigg).
\]
Since we may assume by Example \ref{ex-xleE} that $x<2\,e_{c}(a) \le 2E$, it follows that
\[
x Y<\frac{\epsilon +1}{\delta}\,E.
\]
On the other hand, Lemma \ref{lem-Delta} with $c'=c$ says that $xY-1=\Delta=jE$ for some $j \in \mathbb N$.
These together easily imply the asserted bound for $j$. % by the arbitrary choice of $\epsilon$.
\par
(ii-2) 
Since $c^Z=a+b^Y<2b^Y$, one has
\[
c^{\delta Z}< (2b^Y)^\delta= 2^\delta \cdot b^{\delta Y}.
\]
Since $a<c^{z/x}$, it follows from \eqref{agebdeltaY} that 
\[
c^{\delta Z}< 2^\delta \cdot a<2^\delta \cdot c^{z/x},
\]
so that $c^{\,\delta Z-z/x}< 2^\delta$, leading to
\[
Z <\frac z{\delta x} +\frac{\log 2}{\log c},
\]
whence
%\begin{equation}\label{ineq-z}
\[
\left(1-\frac 1{\delta x}\right)z < \frac{\log 2}{\log c} - (Z-z) \le 1.
\]
%end{equation}
From now on, assume that $1-\frac 1{\delta x}>0$, i.e., $x>1/\delta$.
Then the above displayed inequalities give
\[
z<\frac{1}{1-\frac 1{\delta x}}.
\]
Further, since $x \ll_c 1$, one obtains $z \ll_{c,\delta}1$, so that $\max\{a,b\}<c^z \ll_{c,\delta} 1$.
This completes the proof.
\end{proof}

\begin{ex}\rm
In Proposition \ref{prop-systemX1-bounds}, taking $\delta=0.51$ and $\epsilon=0.01$, we have
\[
\begin{cases}
\,Y <5143 \log c& \text{if $a<b^{0.51Y}\,;$}\\
\,x=1, \ \ Y=1+E & \text{if $a \ge b^{0.51Y}$,}\\
%\,x<5040 \log c, \ \ z=Z<\frac{\log 8}{\log c} & \text{if $a \ge b^{3Y/4}$ and $x \ge 2$,}
\end{cases}
\]
whenever $\max\{a,b\}$ exceeds some constant which depends only on $c$ and is effectively computable.
\end{ex}

%%%%%%%%%%%%%%%%%%%
\subsubsection{Case where $Y=1$}
%%%%%%%%%%%%%%%%%%%

Consider the following system of two equations:
\begin{eqnarray}
&a^x+b=c^z,\label{eff-1st-y1Y1}\\
&a^X+b=c^Z \label{eff-2nd-y1Y1}
\end{eqnarray}
with $x<X$ and $z<Z$.
It is not hard to see that $a^X>b$.
%\[
%a^x<6000b.
%\]

\begin{prop}\label{prop-systemY1-bounds}
For any number $\delta$ with $0<\delta<1,$ the following hold.
\begin{itemize}
\item[\rm (i)] 
If $a^{\delta X}>b,$ then either 
\[
X < \min\biggl\{s_0,\,\frac{2520}{1-\delta}\biggr\} \cdot \log c,
\]
or
\[
s_0 \le s<\frac{25.2}{\,1-\delta\,}\,\big( \log (2s+1) +0.38\,\bigr)^2,
\]
where $s=X/\log c$ and $s_0$ is the same as that in Proposition $\ref{prop-systemX1-bounds}.$
\item[\rm (ii)] 
Assume that $a^{\delta X} \le b.$ 
Let $\epsilon$ be any positive number with $\epsilon<1.$
Then the following hold.
\begin{itemize}
\item[\rm (ii-1)] 
If $\max\{a,b\}$ exceeds some constant which depends only on $c$ and $\epsilon$ and is effectively computable, then
\[
\,X=x+jE
\]
for some positive integer $j$ with $j<(1+\epsilon)/\delta.$
\item[\rm (ii-2)] 
If $\max\{a,b\}$ exceeds some constant which depends only on $c, \delta$ and $\epsilon$ and is effectively computable, then
\[
x < \left(\frac{1}{\delta}-1\right)E+\epsilon.
\]
\end{itemize}
\end{itemize}
\end{prop}

\begin{proof}
(i) The assertion can be proved similarly to the proof of Proposition \ref{prop-systemX1-bounds}\,(i).
\par
(ii-1) Assume that
\begin{equation}\label{adeltaXleb}
a^{\delta X} \le b.
\end{equation}
Since $b<c^z$, and $c^z \ll_{c}a^E$ by Lemma \ref{lem-heggcz}, one has
\[
b \ll_{c} a^E.
\]
This together with \eqref{adeltaXleb} implies that
\[
a^{\delta X-E} \ll_{c} 1.
\]
If $\delta X-E \ge \epsilon$, then $a^{\epsilon} \ll_{c} 1$, so that $a \ll_{c,\epsilon} 1$, leading to $\max\{a,b\}<c^z \ll_c a^E \ll_{c,\epsilon} 1$.
Thus we may assume that $\delta X-E < \epsilon$, i.e.,
\[
X<\frac{\epsilon+E}{\delta}.
\]
On the other hand, $X-x=\Delta=j E$ for some $j \in \mathbb N$.
These together easily imply the asserted bound for $j$.
\par
(ii-2) Since $a^X>b$, it follows from \eqref{adeltaXleb} that 
\[
c^Z=a^X+b<2a^X \le 2 b^{1/\delta}<2c^{z/\delta}.
\]
On the other hand, reducing equations \eqref{eff-1st-y1Y1}, \eqref{eff-2nd-y1Y1} modulo $a^x$ together implies that $c^{Z-z} \equiv 1 \pmod{a^x}$, in particular, 
\[
c^{Z-z}>a^x.
\]
Since $c^z \ll_{c} a^E$, these inequalities together show that
\[
a^x<c^{Z-z}<2\,c^{(1/\delta-1)z} \ll_{c,\delta} a^{(1/\delta-1)E},
\]
so that
\[
a^{x-(1/\delta-1)E} \ll_{c,\delta} 1.
\]
If $x-(1/\delta-1)E \ge \epsilon$, then $a^{\epsilon} \ll_{c,\delta} 1$, so that $a \ll_{c,\delta,\epsilon} 1$, leading to $\max\{a,b\}<c^z \ll_c a^E \ll_{c,\delta,\epsilon} 1$.
This completes the proof.
\end{proof}

\begin{ex}\rm
In Proposition \ref{prop-systemX1-bounds}, taking $\delta=0.51$ and $\epsilon=0.01$, we have
\[
\begin{cases}
\,X <5143 \log c& \text{if $a^{0.51X}>b\,;$}\\
\,x < \frac{49}{51}E + 0.01, \ \ X=x+E & \text{if $a^{0.51X} \le b,$}
\end{cases}
\]
whenever $\max\{a,b\}$ exceeds some constant which depends only on $c$ and is effectively computable.
\end{ex}

%%%%%%%%%%%%%%%%%%
\section{Concluding remarks}
%%%%%%%%%%%%%%%%%%

According to the contents of Sections \ref{sec-6} and \ref{sec-7}, together with that of \cite[Sec.\,7]{MiyPi3}, as a sufficient condition for making Theorem \ref{th-anyc} effective for certain values of $c$, we can raise the following problem, which is closely related to \cite[Proposition 7.1]{MiyPi3}.

\begin{prob}\rm
Let $c \ne 13$ be any prime of the form $c=2^r \cdot 3+1$ with some positive integer $r$. 
For each positive integer $Y$ with $10 \le Y<3520 \log c$ and $Y \equiv 4 \pmod{6}$, and for each positive integer $R$ with $R<Y$, find a positive constant $\epsilon$ and a positive constant $C=C(c,Y,R,\epsilon)$ such that the inequality 
\[
\left| c^{\,R/Y} - \frac{p}{q} \right| > \frac{C}{q^{Y-2-\epsilon}}
\]
holds for any positive integers $p$ and $q$ with the property that $q$ is a power of $c$.
\end{prob}

%\subsection*{Acknowledgements}
%The authors would like to thank Reese Scott for pointing out an error in the proof of Theorem \ref{th-x2y1X1Y2} of an earlier draft.


\begin{thebibliography}{2}
%
\bibitem[Ben]{Ben_cjm_01}
M.A.~Bennett,
{\it On some exponential equations of S. S. Pillai},
Canad.\ J.\ Math. {\bf 53}\,(2001), no. 5, 897--922.
%
\bibitem[BenMihSik]{BenMihSik}
M.A.~Bennett, P. Mih\u{a}ilescu and S. Siksek,
{\it The generalized Fermat equation}, in Springer volume Open Problems in Mathematics, 2016, 173--205.
%
%\bibitem[BeSi]{BeSi}
%M.A. Bennett and S. Siksek,
%{\it Differences between perfect powers$:$ prime power gaps}, Algebra Number Theory {\bf 17}\,(2023), no. 10, 1789--1846.
%
\bibitem[Beu]{Beu}
F. Beukers, 
{\it The Diophantine equation $Ax^p+By^q=Cz^r$}, 
Duke Math. J.\ {\bf 91}\,(1998), no. 1, 61--88. 
%
%\bibitem[BeuSc]{BeSc}
%F. Beukers and H.P. Schlickewei, 
%{\it The equation $x+y=1$ in finitely generated groups}, 
%Acta Arith.\ {\bf 78}\,(1996), 189--199. 
%
\bibitem[BiBu]{BiBu}
Y. Bilu and Y. Bugeaud,
{\it D\'emonstration du th\'eor\'eme de Baker-Feldman via les formes lin\'eaires en deux logarithmes},
J.\ Theor.\ Nombres Bordeaux {\bf 12}\,(2000), no. 1, 13--23.
%
\bibitem[Bu]{Bu}
Y. Bugeaud,
{\it On the greatest prime factor of $ax^m+by^n,$ II.},
Bull.\ London Math.\ Soc.\ {\bf 32}\,(2000), no. 6, 673--678.
%
\bibitem[BuLu]{BuLu}
Y. Bugeaud and F. Luca,
{\it On Pillai's Diophantine equation},
New York J.\ Math.\ {\bf 12}\,(2006), 193--217.
%
\bibitem[EvGy]{EvGy}
J.H.~Evertse and K. Gy{\H o}ry,
{\it Unit Equations in Diophantine Number Theory},
Cambridge University Press, Cambridge, 2015.
%
\bibitem[Fe]{Fe}
N.I.~Fel'dman,
{\it An effective refinement of the exponent in Liouville's theorem}, 
Math.\ USSR, Izv. {\bf 5}\,(1971), no. 5, 985--1002.
%
\bibitem[Ge]{Ge}
A.O.~Gel'fond, 
{\it Sur la divisibilit\'e de la diff\'erence des puissance de deux nombres entiers par une puissance d'un id\'eal premier}, 
Mat. Sb. {\bf 7}\,(1) (1940), 7--25.
%
\bibitem[Gu]{Gu}
R.K.~Guy,
{\it Unsolved Problems in Number Theory},
Springer, 2004.
%
\bibitem[HuLe]{HuLe}
Y.-Z. Hu and M.-H. Le,
{\it A note on ternary purely exponential diophantine equations},
Acta Arith.\ {\bf 171}\,(2015), no. 2, 173--182.
%
\bibitem[La]{La}
M. Laurent,
{\it Linear forms in two logarithms and interpolation determinants II},
Acta Arith.\ {\bf 133}\,(2008), 325--348.
%
\bibitem[Li]{Li}
J. Liouville,
{\it Sur des classes tr\'es-\'etendues de quantit\'es dont la valeur n'est ni alg\'ebrique, ni m\^{e}me r\'eductible \`{a} des irrationnelles alg\'ebriques},
C.\ R.\ Acad.\ Sci.\ Paris.\ {\bf 18}\,(1844), 883--885.
%
\bibitem[Ma]{Ma}
K. Mahler, {\it Zur Approximation algebraischer Zahlen I\,$:$ \"Uber den gr\"ossten Primteiler bin\"arer Formen}, 
Math.\ Ann.\ {\bf 107}\,(1933), 691--730.
%
\bibitem[MiyPi]{MiyPi}
T. Miyazaki and I. Pink,
{\it Number of solutions to a special type of unit equations in two unknowns}, 
Amer.\ J.\ Math.\ {\bf 146}\,(2024), no. 2, 295--369.
%
\bibitem[MiyPi2]{MiyPi2}
---\!---,
{\it Number of solutions to a special type of unit equations in two unknowns, I\hspace{-1.2pt}I}, 
Res.\ Number Theory {\bf 10}\,(2024), no. 2, 41 pp.
%
\bibitem[MiyPi3]{MiyPi3}
---\!---,
{\it Number of solutions to a special type of unit equations in two unknowns, I\hspace{-1.2pt}I\hspace{-1.2pt}I}, preprint 2024, arXiv:2403.20037 (accepted for publication in Math. Proc. Cambridge Philos. Soc.).
%
\bibitem[Ri]{Rid}
D. Ridout,
{\it Rational approximations to algebraic numbers},
Mathematika {\bf 4}\,(1957), 125--131.
%
\bibitem[Sch]{Sc}
W.M.~Schmidt,
{\it Diophantine approximation}, 
Lecture Notes in Mathematics, vol. 785, Springer, 1980.
%
\bibitem[Sch2]{Sc2}
W.M.~Schmidt,
{\it Diophantine approximations and Diophantine equations}, 
Lecture Notes in Mathematics, vol. 1467, Springer, 1991.
%
\bibitem[Sco]{Sc_JNT_93}
R. Scott,
{\it On the equations $p^x-b^y=c$ and $a^x+b^y=c^z$},
J.\ Number Theory {\bf 44}\,(1993), no. 2, 153--165.
%
\bibitem[ScoSt]{ScSt_PMD_2016}
R. Scott and R. Styer,
{\it Number of solutions to $a^x+b^y=c^z$},
Publ.\ Math.\ Debrecen {\bf 88}\,(2016), no. 1-2, 131--138.
%
\bibitem[ShTi]{ShTi}
T.N.~Shorey and R. Tijdeman,
{\it Exponential Diophantine Equations, Cambridge Tracts in Math.},
vol. 87, Cambridge University Press, Cambridge, 1986.
%
\bibitem[Sie]{Sie} 
W. Sierpi\'nski,
{\it On the equation $3^x+4^y=5^z$},
Wiadom.\ Mat.\ {\bf 1}\,(1955/56), 194--195 (in Polish).
%
\bibitem[Wa]{Wa}
M. Waldschmidt, 
{\it Perfect Powers$:$ Pillai's works and their developments}, 
R. Balasubramanian, R. Thangadurai (Eds.), Collected Works of S. Sivasankaranarayana Pillai, Collected Works Series, vol. 1, Ramanujan Mathematical Society, Mysore, 2010.
%
\end{thebibliography}
\end{document}